\documentclass[reqno]{amsart}
\usepackage{amsmath,amssymb}
\usepackage{txfonts}
\usepackage{rrproof}

\usepackage{ifthen}

\newif\ifpdf
\ifx\pdfoutput\undefined
  \pdffalse
\else
  \ifx\pdfoutput\relax
  \else
    \ifnum\pdfoutput>0
      \pdftrue
    \fi
  \fi
\fi

\newif\ifdraft 
\newif\ifshort 

\newif\iffalse 
\falsefalse

\newenvironment{notation*}{\par\noindent\bf Notation:~\rm}{}

\newcounter{mboldlevel}
\let\oldmb\mathbf
 \renewcommand{\mathbf}[1]{
   \addtocounter{mboldlevel}{1}
   \oldmb{#1}
   \addtocounter{mboldlevel}{-1}
}

\newenvironment{leftalign*}{%
\begin{equation*}\begin{array}{l l r}
}{\end{array}\end{equation*}}

\newcommand{\mb}[1]{\ensuremath{\mathbf{#1}}}
\newcommand{\mc}[1]{\ensuremath{\mathcal{#1}}}

\newcommand{\ms}[1]{\lower -0.3ex\hbox{\ensuremath{\scriptstyle #1}}}

\providecommand{\url}[1]{$\langle$\texttt{#1}$\rangle$}
\providecommand{\doi}[1]{\url{doi:#1}}

\newcommand{\sn}[1]{\mb{#1}}

\newcommand{\Gm}{\Gamma}
\newcommand{\Dl}{\Delta}
\newcommand{\Sg}{\Sigma}

\newcommand{\hs}[1]{\mc{#1}}         
\newcommand{\hH}{\hs{H}}                  

\newcommand{\p}{\ifmmode ^{\prime} \else \(^{\prime}\) \fi}
\newcommand{\pp}{\ifmmode ^{\prime\prime}
  \else \(^{\prime\prime}\) \fi}
\newcommand{\ppp}{\ifmmode ^{\prime\prime\prime}
  \else \(^{\prime\prime\prime}\) \fi}

\newcommand{\imp}{\!\supset}                

\newcommand{\seq}{\!\Rightarrow\!}            

\usepackage{slashed}

\newcommand{\HIDE}[1]{}

\usepackage{amsthm}    
\usepackage{makecmds}  

\newcommand{\theoremalias}[2]{
 \makeenvironment{#1}{\begin{#2}}{\end{#2}}
}

  \theoremstyle{plain}
  
  \newtheorem*{thm*}{Theorem}
  \newtheorem*{lem*}{Lemma}
  \newtheorem*{prop*}{Proposition}
  \newtheorem*{cor*}{Corollary}

  \newtheorem*{examp*}{Example}

  \theoremstyle{definition}
  
  \newtheorem*{defn*}{Definition}
  \newtheorem*{notat*}{Notation} 
  \newtheorem*{term*}{Terminology}

  \newtheorem*{conj*}{Conjecture}

  \theoremstyle{remark}
  
  \newtheorem*{rem*}{Remark}
  \newtheorem*{note*}{Note}

\theoremalias{definition}{defn}
\theoremalias{notation}{notat}
\theoremalias{terminology}{term}

\theoremalias{lemma}{lem}
\theoremalias{theorem}{thm}
\theoremalias{corollary}{cor}
\theoremalias{proposition}{prop}
\theoremalias{conjecture}{conj}
\theoremalias{remark}{rem}
\theoremalias{example}{examp}

\theoremalias{definition*}{defn*}
\theoremalias{notation*}{notat*}

\theoremalias{lemma*}{lem*}
\theoremalias{theorem*}{thm*}
\theoremalias{corollary*}{cor*}
\theoremalias{proposition*}{prop*}
\theoremalias{conjecture*}{conj*}
\theoremalias{remark*}{rem*}
\theoremalias{example*}{examp*}

\renewcommand{\imp}{\to}

\newcommand{\poss}{\Diamonddot}
\newcommand{\necc}{\boxdot}

\newcommand{\fw}[1]{\emph{#1}} 

\usepackage{float}    
\floatstyle{boxed}
\restylefloat{figure}
\floatstyle{ruled}
\restylefloat{table}
\title{Modality for Free:
 Notes on Adding the Tarskian M\"oglichkeit to Substructural Logics}
\author{Robert Rothenberg
}

\ifdraft
  \usepackage{datetime}
  \date{*** DRAFT *** \usdate \today\ \currenttime\ *** DRAFT ***}
\else
  \usepackage{datetime}
  \date{\usdate \today}
\fi

\begin{document}

\begin{abstract} We briefly examine the modal formulae that can be
derived in Multiplicative Additive Linear Logic (\sn{MALL}) and some extensions
by using Tarksi's extensional modal operators.  We also breifly
compare this with a substructural form of the modal logic \sn{K}.
\end{abstract}

\maketitle

\section{Introduction}

The Tarskian \fw{m\"oglichkeit} (literally, ``possibility'' in German)
is a modal operator that was  introduced by {\L}ukasiewicz
(and attributed to Tarski) in \cite[\S 7]{Luka1930}.
This modal operator is unusual in that it is an \emph{extensional} one,
defined in terms of other connectives in {\L}ukasiewicz's many-valued logics:
\begin{equation}
 \poss A =_{def} \neg A \imp A\label{eq:diamond_def}
\end{equation}
The modal logic
that results from this definition is unusual, in part because
of the theorems such as:
\begin{equation}\label{hated_theorem}
 (\poss A \land \poss B) \imp \poss(A\land B)
\end{equation}
In the case where $B = \neg A$, theorem \eqref{hated_theorem} appears
to be paradoxical, if not absurd, and largely because of this, the Tarskian
\fw{m\"oglichkeit} has been a footnote in the history of modal logic.
Most of the analyses that we are aware of has been for the 3-valued
logic, in \cite{LewiLang1959}, 
\cite{HughCres1968} (but omitted from \cite{HughCres1996}), \cite{Dumi1977v4},
and \cite{BullSege1986v2}, and it is generally critical.
An application of the $m>3$-valued logics to describing $m$-state
systems was suggested in \cite{Turq1997-BSL}, and an application of the
infinite-valued logic applied to modelling degrees of believability
was suggested by the current author in \cite{Roth2005}.

However, the infinite-valued logic can be seen as an extension of
Affine Logic \cite{CiabLuch1997}, and many of the modal formulae
derivable in the infinite-valued logic are derivable in weaker substructural
logics.  We give an overview of some of the \emph{formal} properties below,
by noting modal rules and formulae in the corresponding logics.
We make no claims about the applications of the Tarskian m\"oglichkeit.

\section{Multiplicative Additive Linear Logic}

The sequent rules for \sn{GMALL}, a calculus for
Multiplicative Additive Linear Logic (\sn{MALL})
\cite{Gira1987} are given in Figure~\ref{fig:MALL}, using notation
similar to \cite{Troe1992}---in particular, we use $\oplus$ for \emph{par}
(multiplicative disjunction) and $\lor$ for \emph{plus} (additive disjunction).

\begin{figure}[ht]

$$
\infer[L0]{0 \seq }{}
\qquad
\infer[R0]{\Gm\seq\Dl,0}{
 \Gm\seq\Dl
}
\qquad
\infer[L1]{1,\Gm\seq\Dl}{
 \Gm\seq\Dl
}
\qquad
\infer[R1]{\seq 1}{}
$$

$$
\infer[L\bot]{\bot,\Gm\seq\Dl}{}
\qquad
\infer[R\top]{\Gm\seq\Dl,\top}{}
$$

$$
\infer[L\otimes]{A\otimes B,\Gm\seq\Dl}{
 \Gm,A,B\seq\Dl
}
\qquad
\infer[R\otimes]{\Gm,\Gm'\seq\Dl,\Dl',A\otimes B}{
 \Gm\seq\Dl,A &
 \Gm'\seq\Dl',B
}
$$

$$
\infer[L\oplus]{A\oplus B,\Gm,\Gm'\seq\Dl,\Dl'}{
 A,\Gm\seq\Dl &
 B,\Gm'\seq\Dl'
}
\qquad
\infer[R\oplus]{\Gm\seq\Dl,A\oplus B}{
 \Gm\seq\Dl,A,B
}
$$

$$
\infer[L\land_1]{A\ \land\ B,\Gm\seq\Dl}{
 A,\Gm\seq\Dl
}
\quad
\infer[L\land_2]{A\ \land\ B,\Gm\seq\Dl}{
 B,\Gm\seq\Dl
}
\quad
\infer[R\land]{\Gm\seq\Dl,A\ \land\ B}{
 \Gm\seq\Dl,A & \Gm\seq\Dl,B
}
$$

$$
\infer[L\lor ]{A\lor B,\Gm\seq\Dl}{
 A,\Gm\seq\Dl & B,\Gm\seq\Dl
}
\quad
\infer[R\lor _1]{\Gm\seq\Dl,A\lor B}{
 \Gm\seq\Dl,A
}
\quad
\infer[R\lor _2]{\Gm\seq\Dl,A\lor B}{
 \Gm\seq\Dl,B
}
$$

$$
\infer[L\neg]{\neg A,\Gm\seq\Dl}{
 \Gm\seq\Dl,A
}
\qquad
\infer[R\neg]{\Gm\seq\Dl,\neg A}{
 A,\Gm\seq\Dl
}
$$

$$
\infer[L\imp]{A\imp B,\Gm,\Gm'\seq\Dl,\Dl'}{
 \Gm\seq\Dl,A &
 B,\Gm'\seq\Dl'
}
\qquad
\infer[R\imp]{\Gm\seq\Dl,A\imp B}{
 A,\Gm\seq\Dl,B
}
$$

\vspace{1ex}
Note that the rules for $\imp$ can be derived using the definition
$A\imp B =_{def} \neg A\oplus B$.
(Rules for additive implication are omitted but can be derived similarly.)

\caption{Rules for \sn{GMALL}.}\label{fig:MALL}
\end{figure}

\begin{figure}[ht]
$$
\vcenter{
 \infer[L\poss]{\poss A,\Gm,\Gm'\seq\Dl,\Dl'}{
  A,\Gm\seq\Dl &
  A,\Gm'\seq\Dl'
 }
}
\qquad
\vcenter{
 \infer[R\poss]{\Gm\seq\Dl,\poss A}{
  \Gm\seq\Dl,A,A
 }
}
$$

$$
\vcenter{
 \infer[L\necc]{\necc A,\Gm\seq\Dl}{
  A,A,\Gm\seq\Dl
 }
}
\qquad
\vcenter{
 \infer[R\necc]{\Gm,\Gm'\seq\Dl,\Dl',\necc A}{
  \Gm\seq\Dl,A &
  \Gm'\seq\Dl',A
 }
}
$$
\caption{Derived rules for Tarskian modalities in \sn{GMALL}.}\label{fig:modal_rules}
\end{figure}

The modal rules (Figure~\ref{fig:modal_rules})
are derived in a straightforward manner.
(The corresponding box operator is defined as the dual of diamond operator
$\necc A =_{def} \neg\poss\neg A$.)

\begin{remark}
  Because of the symmetries that occur in many of the proofs given in
  this paper, the following non-branching forms of the modal rules
  will be used for brevity:
$$
\vcenter{
 \infer[L\poss]{\poss A,\Gm,\Gm\seq\Dl,\Dl}{
  A,\Gm\seq\Dl
 }
}
\qquad
\vcenter{
 \infer[R\necc]{\Gm,\Gm\seq\Dl,\Dl,\necc A}{
  \Gm\seq\Dl,A
 }
}
$$
\end{remark}

\begin{proposition}\label{prop:alt_defs}
The following equivalences hold in \sn{MALL}:
\begin{eqnarray}
 \poss A \equiv A\oplus A  \\
 \necc A \equiv A\otimes A
\end{eqnarray}
\begin{proof} Straightforward.
\end{proof}
\end{proposition}

\begin{remark}
The equivalences in Proposition~\ref{prop:alt_defs}
may be used as alternative definitions of the Tarskian modalities.
\end{remark}


\begin{proposition} The following are derivable in \sn{GMALL}:
\begin{eqnarray}
 \neg\necc\bot \label{eq:AntiIntD} \\
 \poss\top \label{eq:IntD} \\
 \neg\necc(A\land\neg A) \label{eq:AntiDN} \\
 \poss(A\lor\neg A) \label{eq:DN}
\end{eqnarray}
where \eqref{eq:IntD} corresponds to the intuitionistic modal axiom
\rn{D} \cite{Simp1994-PhD}.
\begin{proof} Straightforward.
\end{proof}
\end{proposition}

\begin{proposition}
The \rn{K\necc} rule and its dual
\begin{equation*}\label{eq:N}
\vcenter{
\infer[K\necc]{\necc\Gm\seq\necc A}{
 \Gm\seq A
}}
\qquad
\vcenter{
\infer[K\poss]{\poss A\seq\poss\Dl}{
 A\seq\Dl
}}
\end{equation*}
are derivable in \sn{GMALL}.
\begin{proof} Straightforward.
\end{proof}
\end{proposition}

\begin{proposition}[Distribution Theorems]
The following are derivable in \sn{GMALL}:
\begin{eqnarray}
 \necc(A\imp B)\imp (\necc A\imp\necc B) \label{eq:K} \\
 \necc(A\land  B)\imp(\necc A \land \necc B) \label{eq:D1} \\
 \poss(A\land  B)\imp(\poss A \land \poss B)\label{eq:D2} \\
 (\necc A\lor\necc B)\imp\necc(A\lor B) \\
 (\poss A\lor \poss B)\imp\poss(A\lor B) \label{eq:D4}
\end{eqnarray}
where \eqref{eq:K} corresponds to the \rn{K} axiom.
\begin{proof} Straightforward.
\end{proof}
\end{proposition}


\section{A Comparison of \sn{MALL} with Substructural-\sn{K} (\sn{KMALL})}

\begin{definition}[Substructural-\sn{K}]
Let Substructural-\sn{K} (\sn{KMALL}) be \sn{MALL} augmented by extending the
language of formulae with $\Box A$. We obtain a calculus \sn{GKMALL} for
\sn{KMALL} by adding the following rule to \sn{GMALL}:
$$
\infer[K\Box]{\Box\Gm\seq\Box A}{
 \Gm\seq A
}
$$
which corresponds to adding to \sn{MALL} the necessity rule
$$
\infer[N]{\vdash \Box A}{
 \vdash A
}
$$
and the \rn{K} axiom, $\Box(A\imp B)\imp(\Box A\imp\Box B)$.
\end{definition}

\begin{remark}
Cut elimination for \sn{GKMALL} is shown in Appendix~\ref{sec:KMALL}.
\end{remark}

\begin{theorem}\label{eq:KMALL_MALL}
If~ $\sn{GKMALL}\vdash \Gm\seq\Dl$, then
$\sn{GMALL}\vdash[\necc/\Box]\Gm\seq [\necc/\Box]\Dl$.
\begin{proof} By induction on the derivation height. Note that
$\Box$-formulae are be introduced into a \sn{KMALL} derivation
either by axioms \rn{L\bot} and \rn{R\top}, or by the \rn{K\Box} rule.
Instances of \rn{K\Box} are replaced by instances of \rn{K\necc}.
\end{proof}
\end{theorem}

\begin{remark} $\sn{GKMALL}\nvdash \neg\Box\bot, \Diamond\top, \neg\Box(A\land \neg A)$ and $\Diamond(A\lor\neg A)$.  A form of the \rn{K\Box} rule that
allows for empty succedents, e.g.
$$
\infer[K\Box']{\Box\Gm\seq\Box\Dl}{
 \Gm\seq\Dl
}
$$
where $|\Dl|\le 1$,
would allow for the derivation of $\neg\Box\bot$ \eqref{eq:AntiIntD} and
$\Diamond\top$ \eqref{eq:IntD}
but not $\neg\Box(A\land \neg A)$ \eqref{eq:AntiDN} and
$\Diamond(A\lor\neg A)$ \eqref{eq:DN}.
\end{remark}

\begin{remark}
The converse of \eqref{eq:D1},
$(\necc A \land\necc B)\imp\necc(A\land B)$ is {not}
derivable in either \sn{GMALL} or \sn{GKMALL}.
\end{remark}

\section{Multiplicative Additive Linear Logic with Mingle}

Linear Logics with mingle are discussed in
\cite{DBLP:journals/jolli/Kamide02,DBLP:journals/logcom/Kamide05}.
The mingle (also called ``merge'' or ``mix'') rule is:
$$
\vcenter{
 \infer[M]{\Gm,\Gm'\seq\Dl,\Dl'}{
  \Gm\seq\Dl &
  \Gm'\seq\Dl'
 }
}
$$
\begin{proposition}\label{prop:MALL+M}
The following are derivable in \sn{GMALL} + \rn{M}:
\begin{eqnarray}
 \necc A\imp\poss A \label{eq:D} \\
 (\poss A\imp\poss B)\imp\poss(A\imp B) \label{eq:K'}
\end{eqnarray}
where \eqref{eq:D} corresponds to a form of the \rn{D} axiom.
\end{proposition}
\begin{remark}
  We note that the formulae \eqref{eq:D} and \eqref{eq:K'} can be
  derived using anti-contraction (duplication) rules as well:
$$
\vcenter{
 \infer[LC^{-1}]{A,A,\Gm\seq\Dl}{
  A,\Gm\seq\Dl
 }
}
\qquad
\vcenter{
 \infer[RC^{-1}]{\Gm\seq\Dl,A,A}{
  \Gm\seq\Dl,A
 }
}
$$
\end{remark}

\section{Affine Logic}

Affine Logic \cite{Gris1982}, also called Affine Multiplicative Additive Linear
Logic (\sn{AMALL}), is \sn{MALL} augmented with the weakening axiom,
$(A\imp 1)\land(0\imp A)$. The
corresponding calculus, \sn{GAMALL} is obtained by adding weakening rules
to \sn{GMALL}:
$$
\infer[LW]{A,\Gm\seq\Dl}{
 \Gm\seq\Dl
}
\qquad
\infer[RW]{\Gm\seq\Dl,A}{
 \Gm\seq\Dl
}
$$
\begin{remark}
In Affine Logic, $0 \equiv \bot$ and $1 \equiv \top$.
\end{remark}

\begin{proposition} If $\sn{GMALL}+\rn{M}\vdash A$, then $\sn{GAMALL}\vdash A$.
\begin{proof} \rn{M} is admissible in \sn{GAMALL}.
\end{proof}
\end{proposition}

\begin{remark} Hence, the formulae in Proposition~\ref{prop:MALL+M}
are derivable in \sn{GAMALL}.
\end{remark}

\begin{proposition}[Additive Modal Rules]
The following rules can be derived in \sn{GAMALL}:
$$
\vcenter{
 \infer[R\poss']{\Gm\seq\Dl,\poss A}{
    \Gm\seq\Dl,A
 }
}
\qquad
\vcenter{
 \infer[L\boxdot']{\Gm,\boxdot A\seq\Dl}{
  \Gm,A\seq\Dl
 }
}
$$
\begin{proof} Straightforward.
\end{proof}
\end{proposition}

\begin{proposition} The following are derivable in \sn{GAMALL}:
\begin{eqnarray}
 \necc A\imp A \label{eq:T} \\
 A\imp\poss A
\end{eqnarray}
where \eqref{eq:T} corresponds to the \rn{T} axiom.
\begin{proof} Straightforward.
\end{proof}
\end{proposition}

\begin{proposition} The following rules are derivable in \sn{GAMALL}:
$$
\infer[R\poss R\necc]{\seq\poss(A\imp\necc B)}{
 A\seq B
}
\qquad
\infer[R\poss L\poss]{\seq\poss(\poss A\imp B)}{
 A\seq B
}
$$
\begin{proof} Straightforward.
\end{proof}
\end{proposition}

\begin{proposition}
The following are derivable in \sn{GAMALL}:
\begin{eqnarray}
 \poss(\necc A\imp\necc\necc A) \label{eq:S4} \\
 \poss(\poss A\imp\necc\poss A) \label{eq:S5} \\
 \poss(A\imp\necc\poss A) \label{eq:B} \\
 \poss(\poss A\imp A) \\
 \poss(A\imp\necc A) \\
 \poss(A\imp\necc B)\imp\poss(\necc A\imp B)
\end{eqnarray}
where \eqref{eq:S4}, \eqref{eq:S5} and \eqref{eq:B} are ``$\poss$-forms'' of
the \rn{S4}, \rn{S5} and \rn{B} axioms, respectively.
\begin{proof} Straightforward.
\end{proof}
\end{proposition}

\section{Strict Logic}

A calculus \sn{GSLL} for Strict Linear Logic (\sn{SLL}) is obtained by
adding to \sn{GMALL} the contraction rules:
$$
\infer[LC]{A,\Gm\seq\Dl}{
 A,A,\Gm\seq\Dl
}
\qquad
\infer[RC]{\Gm\seq\Dl,A}{
 \Gm\seq\Dl,A,A
}
$$

\begin{proposition}
The following can be derived in \sn{GSLL}:
\begin{eqnarray}
 A\imp\necc A \\
 \poss A\imp A \\
 \necc A\imp\necc \necc A \label{eq:S4-SLL} \\
 \poss\poss A\imp A\poss \\
 \poss A\imp\necc\poss A \label{eq:B-SLL} \\
 A\imp\necc\poss A \label{eq:S5-SLL}
\end{eqnarray}
where \eqref{eq:S4-SLL}, \eqref{eq:B-SLL} and \eqref{eq:S5-SLL}
correspond to the \rn{S4}, \rn{B} and \sn{S5} axioms, respectively.
\begin{proof} Straightforward.
\end{proof}
\end{proposition}

\begin{proposition}
The following can be derived in \sn{GSLL}:
\begin{eqnarray}
(\necc A\land \necc B)\imp\necc(A\land B) \label{eq:D1'} \\
{(\poss A\land \poss B)}{\imp\poss(A\land B)} \label{eq:D2'} \\
{\necc(A\lor B)}{\imp(\necc A\lor \necc B)}\\
\poss(A\lor B)\imp(\poss A\lor \poss B) \label{eq:D4'}
\end{eqnarray}
\begin{proof} Straightforward.
\end{proof}
\end{proposition}

\begin{remark} These are the converse of formulae \eqref{eq:D1'}
through \eqref{eq:D4'}.
Note that \eqref{eq:D2'} is the same formulae as \eqref{hated_theorem}
mentioned in the introduction. Indeed \eqref{eq:D1'} may also
be considered paradoxical.
\end{remark}

\section{Involutive Uninorm Logic}

Involutive Uninorm Logic (\sn{IUL}) \cite{MetcMont2007-JSL72(3)}
is a substructural fuzzy logic, and
has a hypersequent calculus \sn{GIUL} (Figure~\ref{fig:GIUL})
based on a hyperextension of \sn{GMALL} \cite{Avro1991}
and the communication (\rn{Com}) rule.

\begin{figure}[ht]

$$
\infer[L0]{0 \seq }{}
\qquad
\infer[R0]{\hH~|~\Gm\seq\Dl,0}{
 \hH~|~\Gm\seq\Dl
}
\qquad
\infer[L1]{\hH~|~1,\Gm\seq\Dl}{
 \hH~|~\Gm\seq\Dl
}
\qquad
\infer[R1]{\seq 1}{}
$$

$$
\infer[L\bot]{\bot,\Gm\seq\Dl}{}
\qquad
\infer[R\top]{\Gm\seq\Dl,\top}{}
$$

$$
\infer[L\otimes]{\hH~|~A\otimes B,\Gm\seq\Dl}{
 \hH~|~\Gm,A,B\seq\Dl
}
\qquad
\infer[R\otimes]{\hH~|~\Gm,\Gm'\seq\Dl,\Dl',A\otimes B}{
 \hH~|~\Gm\seq\Dl,A &
 \hH~|~\Gm'\seq\Dl',B
}
$$

$$
\infer[L\oplus]{\hH~|~A\oplus B,\Gm,\Gm'\seq\Dl,\Dl'}{
 \hH~|~A,\Gm\seq\Dl &
 \hH~|~B,\Gm'\seq\Dl'
}
\qquad
\infer[R\oplus]{\hH~|~\Gm\seq\Dl,A\oplus B}{
 \hH~|~\Gm\seq\Dl,A,B
}
$$

$$
\infer[L\land_1]{\hH~|~A\ \land\ B,\Gm\seq\Dl}{
 \hH~|~A,\Gm\seq\Dl
}
\quad
\infer[L\land_2]{\hH~|~A\ \land\ B,\Gm\seq\Dl}{
 \hH~|~B,\Gm\seq\Dl
}
\quad
\infer[R\land]{\hH~|~\Gm\seq\Dl,A\ \land\ B}{
 \hH~|~\Gm\seq\Dl,A & \hH~|~\Gm\seq\Dl,B
}
$$

$$
\infer[L\lor ]{\hH~|~A\lor B,\Gm\seq\Dl}{
 \hH~|~A,\Gm\seq\Dl & \hH~|~B,\Gm\seq\Dl
}
\quad
\infer[R\lor _1]{\hH~|~\Gm\seq\Dl,A\lor B}{
 \hH~|~\Gm\seq\Dl,A
}
\quad
\infer[R\lor _2]{\hH~|~\Gm\seq\Dl,A\lor B}{
 \hH~|~\Gm\seq\Dl,B
}
$$

$$
\infer[L\neg]{\hH~|~\neg A,\Gm\seq\Dl}{
 \hH~|~\Gm\seq\Dl,A
}
\qquad
\infer[R\neg]{\hH~|~\Gm\seq\Dl,\neg A}{
 \hH~|~A,\Gm\seq\Dl
}
$$

$$
\infer[L\imp]{\hH~|~A\imp B,\Gm,\Gm'\seq\Dl,\Dl'}{
 \hH~|~\Gm\seq\Dl,A &
 \hH~|~B,\Gm'\seq\Dl'
}
\qquad
\infer[R\imp]{\hH~|~\Gm\seq\Dl,A\imp B}{
 \hH~|~A,\Gm\seq\Dl,B
}
$$

$$
\infer[EW]{\hH~|~\Gm\seq\Dl}{
 \hH
}
\qquad
\infer[EC]{\hH~|~\Gm\seq\Dl}{
 \hH~|~\Gm\seq\Dl~|~\Gm\seq\Dl
}
$$

$$
\infer[Com]{\hH~|~\Gm_1,\Gm_2\seq\Dl_1,\Dl_2~|~\Pi_1,\Pi_2\seq\Sg_1,\Sg_2}{
 \hH~|~\Gm_1,\Pi_1\seq\Sg_1,\Dl_1 &
 \hH~|~\Gm_2,\Pi_2\seq\Sg_2,\Dl_2
}
$$


\caption{Rules for \sn{GIUL}.}\label{fig:GIUL}
\end{figure}

\begin{proposition} The following rules are derivable in \sn{GIUL}:\label{prop:additive_contraction}
$$
\infer[L\land']{\hH~|~A\ \land\ B,\Gm\seq\Dl}{
 \hH~|~A,\Gm\seq\Dl~|~B,\Gm\seq\Dl
}
\qquad
\infer[R\lor']{\hH~|~\Gm\seq\Dl,A\lor B}{
 \hH~|~\Gm\seq\Dl,A~|~\Gm\seq\Dl,B
}
$$
\begin{proof} Straightforward, using \rn{EC}.
\end{proof}
\end{proposition}

\begin{samepage}
\begin{proposition} Formulae \eqref{eq:D1'} through
\eqref{eq:D4'} are derivable in \sn{GIUL}.
\begin{proof} Straightforward, using rules from Proposition~\ref{prop:additive_contraction} and \rn{Com}. A proof of \eqref{eq:D1'}:
\begin{equation*}
\vcenter{
 \infer[R\imp]{\seq (\necc A\land \necc B)\imp\necc(A\land B)}{
  \infer[L\land']{\necc A\land \necc B\seq\necc(A\land B)}{
   \infer[L\necc^2]{\necc A\seq\necc(A\land B)~|~\necc B\seq\necc(A\land B)}{
    \infer[R\necc]{A,A\seq\necc(A\land B)~|~B,B\seq\necc(A\land B)}{
     \infer[R\necc]{A\seq A\land B~|~B,B\seq\necc(A\land B)}{
      \infer[R\land]{A\seq A\land B~|~B\seq A\land B}{
       \infer[EW]{A\seq A~|~B\seq A\land B}{
        A\seq A
       }
       &
       \infer[R\land]{A\seq B~|~B\seq A\land B}{
        \infer[Com]{A\seq B~|~B\seq A}{
         A\seq A
         &
         B\seq B
        }
        &
        \infer[EW]{A\seq B~|~B\seq B}{
         B\seq B
        }
       }
      }
     }
    }
   }
  }
 }
}
\end{equation*}
Proofs of  \eqref{eq:D2'} through \eqref{eq:D4'} are similar.



\end{proof}
\end{proposition}
\end{samepage}















\section{Discussion and Future Work}

Much of the content in this paper is straightforward.  However,
the formal properties of the Tarskian m\"oglichkeit are of interest.

Theorem~\ref{eq:KMALL_MALL} is noteworthy, in that all of the
derivable modal sequents derivable in \sn{GKMALL} correspond to
derivable modal sequents in \sn{GMALL} using the Tarskian
m\"oglichkeit.


A semantic characterisation of the Tarskian modalities with respect
to various logics is in process.

A deeper comparison of substructural logics with Tarskian modalities
and their counterpart extensions to \sn{KMALL} is an area of
future investigation.

\section*{Acknowledgements}

We'd like to thank those attending a talk at an LFCS Lab Lunch
about this topic for their comments and suggestions.

\appendix

\section{Substructural-\sn{K} (\sn{KMALL})}
\label{sec:KMALL}

\begin{lemma}[Cut Admissibility] \sn{GKMALL} admits cut
$$
\infer[cut]{\Gm,\Gm'\seq\Dl,\Dl'}{
 \Gm\seq\Dl,A &
 A,\Gm'\seq\Dl'
}
$$
\begin{proof} Note that \sn{GMALL} admits cut \cite{Troe1992}. Adding \rn{K\Box}
to \sn{GKMALL} also admits cut, by induction on the derivation height,
with the following cases:
\begin{enumerate}
\item If the cut formula is not of the form $\Box A$, then permute cut upwards.
\item If the cut formula is not the principal formula on either premiss,
permute the cut upward on that premiss.
\item If the cut formula is the principal formula of either an instance
of \rn{L\bot} or \rn{R\top}, then so is the conclusion of the cut.
\item The remaining case is that both premisses of the cut are the conclusions
of instances of \rn{K\Box}.  The cut is then permuted to the premisses of both
\rn{K\Box} instances, and \rn{K\Box} is applied to the conclusion of the cut.
\end{enumerate}
\end{proof}
\end{lemma}

\bibliographystyle{plain}
\bibliography{macros,logic}

\end{document}